  \documentclass[11pt]{article}
  \usepackage[top=1in,bottom=1in,left=1.5in,right=1.5in]{geometry}
  \usepackage{amsmath,amssymb}
  \usepackage{latexsym}% defines $\Box$ for LaTeX2e
  \usepackage[dvips]{graphicx}
 \newcommand{\submodularf}{w}%SG command for the generic submodular function

   \newcommand{\R}{\mathbb{R}}

   \newcommand{\rH}{\mathrm{H}}

   \newcommand{\hs}{\hspace*{\parindent}}

   \newcommand{\qed}{\hspace*{\fill} $\Box$\\}

   \newtheorem{theo}{\bfseries \hs Theorem}%[section]

   \newtheorem{corol}[theo]{\bfseries \hs Corollary}

   \numberwithin{equation}{section} % Automatically number equations within sections

\renewcommand{\ge}{\geqslant}
\renewcommand{\le}{\leqslant}

\begin{document}

  \title{Nonnegative definite hermitian matrices\\with increasing principal minors}
  \author{
  Shmuel Friedland\footnotemark[1]
  }
  \renewcommand{\thefootnote}{\fnsymbol{footnote}}
  \footnotetext[1]{
  Department of Mathematics, Statistics and Computer Science,
  University of Illinois at Chicago, Chicago, Illinois 60607-7045,
  USA, \texttt{friedlan@uic.edu}. This work was supported by NSF grant DMS-1216393.
  }
  
  \renewcommand{\thefootnote}{\arabic{footnote}}
  \date{January 20, 2013}
  \maketitle
  \begin{abstract}
  A nonnegative definite hermitian  $m\times m$ matrix $A\ne 0$ has increasing principal minors if
  $\det A[I]\le \det[J]$ for $I\subset J$, where $\det A[I]$ is the principal minor of
  $A$ based on rows and columns in the set $I\subseteq \{1,\ldots,m\}$.
  For $m>1$ we show $A$ has increasing principal minors if and only if $A^{-1}$ exists and its diagonal
  entries are less or equal to $1$.
  \end{abstract}



  \section{Introduction and statement of the main result}\label{intro}

  Let $m$ be a positive integer and denote $[m]:=\{1,\ldots,m\}$.
 A real valued function $\submodularf:2^{[m]}\to \R$ defined on all subsets of $[m]$ is called \emph{nondecreasing} if $\submodularf(I)\le \submodularf(J)$ when $I\subset J\subset [m]$.
 It is \emph{submodular} if
 \[ \submodularf(I)+\submodularf(J)\ge \submodularf(I\cup J)+\submodularf(I\cap J)
 \]
 for any two subsets $I,J$ of $[m]$.
 The importance of submodular functions in combinatorial optimization is well known. Several polynomial time algorithms to minimize a submodular function
 under a matroid constraint are known, we refer the reader to the
 survey~\cite{iwata} for more information. The maximization
 of a submodular function under a matroid constraint, and specially,
 under a cardinality constraint,
 $\nu_k(\submodularf):=\max_{I\subset [m], |I|\le k} \submodularf(I)$,
 is also of great interest. 
  For some submodular functions $\submodularf$ the latter problem is NP-hard.
 % A greedy algorithm is often used to approximate $\nu_k(f)$ by $\nu_k^G(f)$.
 However, a classical result~\cite{NWF78} shows
 that when $\submodularf$ is nondecreasing and submodular, the greedy algorithm
 allows one to compute an approximation $\nu_k^G(\submodularf)$ of $\nu_k(\submodularf)$ which is
 such that $\nu_k^G(\submodularf)\ge (1-e^{-1})\nu_k(\submodularf)$.

 Denote by $\rH_{m,+}\supset \rH_{m,++}$ the cone of $m\times m$ nonnegative definite hermitian matrices and its 
 interior consisting of positive definite hermitian matrices respectively. For $I\subseteq [m]$ denote by $A[I]$ the principal submatrix of $A$,
 obtained from $A$ by deleting the rows and columns in the set $[m]\setminus I$.
 Recall that the principal minors of a nonnegative definite matrix satisfy the multiplicative submodularity property:
 \begin{equation}\label{multsubprop}
 \det A[I\cup J]\det A[I\cap J]\le \det A[I]\det A[J], \textrm{ where } I,J\subseteq [m],\;
 A\in\rH_{m,+}.
 \end{equation}
 We assume here that $\det A[\emptyset]=1$.
 In other words, the function $\log(\cdot,A): 2^{[m]}\to \R$ given by
 \begin{equation}\label{deffunc}
 \log (I,A):=\log\det A[I], \quad I\subseteq [m],\; A\in\rH_{m,+}
 \end{equation}
 is submodular.
 This inequality has arisen in the work of several authors.
 It goes back to Gantmacher and Kre{\u\i}n~\cite{gantmacherkrein}
 and Kotelyanski{\u\i}~\cite{kotelyanskii}, see the discussion by Ky Fan~\cite{KyFan67,KyFan68}.
 The classical Hadamard-Fischer inequality for the principal minors of nonnegative definite matrices is obtained when $I\cap J=\emptyset$.
 It is well known that the inequality \eqref{multsubprop} hold also for $M$-matrices, e.g. \cite{Car67}.
 
 \cite[\S5]{FG12}  discusses the CUR approximation \cite{GTZ97} of nonnegative definite hermitian matrix.
 The main problem there is to find a good approximation to the maximum of $\det A[I]$ on all subsets $I$ of $[m]$ of cardinality  $k$.
 Assuming that $A$ has increasing principal minors the greedy algorithm is applied to give an estimate for the $CUR$ approximation.
 It is shown in \cite{FG12} that if all eigenvalues of $A$ are greater or equal $1$ then $A$ has increasing principal minors.
 The purpose of this note is the following theorem.
 \begin{theo}\label{mainthm}  Let $A\in\rH_{m,+}\setminus\{0\}$.  Assume that $m>1$.  Then $A$ has increasing principal minors
 if and only if $A$ is positive definite and all diagonal entries of $A^{-1}$ are less or equal to $1$.
 \end{theo}

 \section{Proof of Theorem \ref{mainthm}}
 Assume first that $A$ has increasing principal minors.
 Suppose that $\det A=0$.  Since $A$ is nonnegative definite we have that $0\le \det A[I]\le \det A=0$ for any nontrivial subset $I$ of $[m]$.
 Hence $A=0$ contrary to our assumption.  Therefore $A$ is positive definite.  Let $B=[b_{ij}]:=A^{-1}$.  Clearly $B$ is positive definite.
 As $A$ has increasing principal minors we deduce that $\det A[[m]\setminus\{i\}]\le \det A$.   Hence $b_{ii}=\frac{\det A[[m]\setminus\{i\}]}{\det A}\le 1$
 for each $i\in [m]$.
 
 It is left to show that if $A$ is positive definite, $B=[b_{ij}]:=A^{-1}$ and $b_{ii}\le 1$ for $i\in [m]$ then $A$ has increasing principal minors.
 We first observe that $B$ has decreasing principal minors, i.e $\det B[I]\ge \det B[J]$ if $I\subset J$.  Indeed, it is enough to consider the case where
 $J=I\cup\{j\}$, where $j\notin I$.  Then the Hadamard-Fischer inequality yields $\det B[J]\le b_{jj} \det B[I] \le \det B[I]$.  Recall the
 Sylvester determinant identity: $\det B[[m]\setminus I]=\frac{\det A[I]}{\det A}$.  Since $B$ has decreasing principal minors it follows that $A$
 has increasing principal minors.  \qed
 
 \begin{corol}\label{tAincrprop}  Let $A\in\rH_{m,++}$ and $m>1$.  Denote $B=[b_{ij}]:=A^{-1}$.  Then $tA, t>0$ has increasing principal minors  
 if and only if $t\ge \max_{i\in [m]} b_{ii}$.
 \end{corol}

 \end{document}